\documentclass[12pt]{article}

\usepackage{amscd,amsmath, amssymb, fancyhdr, epsfig,url}
\usepackage[matrix,arrow]{xy}
\usepackage{diagrams}
\numberwithin{equation}{section}

\newcommand{\version}{version 1.2,\ \   September 10, 2014}

\def\eqref#1{(\ref{#1})}
\newcommand{\goth}{\mathfrak}

\newcommand{\arrow}{{\:\longrightarrow\:}}
\newcommand{\Z}{{\Bbb Z}}
\newcommand{\C}{{\Bbb C}}

\newcommand{\R}{{\Bbb R}}

\renewcommand{\H}{{\Bbb H}}
\newcommand{\6}{\partial}
\def\1{\sqrt{-1}\:}
\newcommand{\restrict}[1]{{\left|_{{\phantom{|}\!\!}_{#1}}\right.}}
\newcommand{\cntrct}                
{\hspace{2pt}\raisebox{1pt}{\text{$\lrcorner$}}\hspace{2pt}}

\makeatletter
\def\x@arrow{\DOTSB\Relbar}
\def\xlongequalsignfill@{\arrowfill@\x@arrow\Relbar\x@arrow}
\newcommand{\xlongequal}[2][]{%
        \ext@arrow 0099\xlongequalsignfill@{#1}{#2}}
\def\xlongrightarrowfill@{\arrowfill@\relbar\relbar\longrightarrow}
\newcommand{\xlongrightarrow}[2][]{%
        \ext@arrow 0099\xlongrightarrowfill@{#1}{#2}}
\makeatother

\newcommand{\calo}{{\cal O}}

\renewcommand{\bar}{\overline}
\renewcommand{\phi}{\varphi}
\renewcommand{\epsilon}{\varepsilon}
\renewcommand{\geq}{\geqslant}
\renewcommand{\leq}{\leqslant}

\newcommand{\im}{\operatorname{im}}

\newcommand{\Av}{\operatorname{\sf Av}}

\newcommand{\Id}{\operatorname{Id}}

\newcommand{\const}{\operatorname{\text{\sf const}}}

\newcommand{\Hol}{\operatorname{Hol}}

\newcommand{\Sym}{\operatorname{Sym}}

\newcommand{\codim}{\operatorname{codim}}

\newcommand{\coker}{\operatorname{coker}}

\newcounter{Mycounter}[section]
\newcounter{lemma}[section]
\renewcommand{\thelemma}{{Lemma \thesection.\arabic{lemma}}}
\newcommand{\lemma}{%
    \setcounter{lemma}{\value{Mycounter}}
    \refstepcounter{lemma}
    \stepcounter{Mycounter}
    {\noindent \bf \thelemma:\ }}

\newcounter{claim}[section]
\renewcommand{\theclaim}{{Claim \thesection.\arabic{claim}}}
\newcommand{\claim}{%
    \setcounter{claim}{\value{Mycounter}}
    \refstepcounter{claim}
    \stepcounter{Mycounter}
    {\noindent \bf \theclaim:\ }}

\newcounter{sublemma}[section]
\newcounter{corollary}[section]
\renewcommand{\thecorollary}{{Corollary \thesection.\arabic{corollary}}}
\newcommand{\corollary}{%
    \setcounter{corollary}{\value{Mycounter}}
    \refstepcounter{corollary}
    \stepcounter{Mycounter}
    {\noindent \bf \thecorollary:\ }}

\newcounter{theorem}[section]
\renewcommand{\thetheorem}{{Theorem \thesection.\arabic{theorem}}}
\newcommand{\theorem}{%
    \setcounter{theorem}{\value{Mycounter}}
    \refstepcounter{theorem}
    \stepcounter{Mycounter}
    {\noindent \bf \thetheorem:\ }}

\newcounter{conjecture}[section]
\newcounter{proposition}[section]
\renewcommand{\theproposition}
      {{Proposition \thesection.\arabic{proposition}}}
\newcommand{\proposition}{%
    \setcounter{proposition}{\value{Mycounter}}
    \refstepcounter{proposition}
    \stepcounter{Mycounter}
    {\noindent \bf \theproposition:\ }}

\newcounter{definition}[section]
\renewcommand{\thedefinition}
      {{Definition~\thesection.\arabic{definition}}}
\newcommand{\definition}{%
    \setcounter{definition}{\value{Mycounter}}
    \refstepcounter{definition}
    \stepcounter{Mycounter}
    {\noindent \bf \thedefinition:\ }}

\newcounter{example}[section]
\renewcommand{\theexample}{{Example \thesection.\arabic{example}}}
\newcommand{\example}{%
    \setcounter{example}{\value{Mycounter}}
    \refstepcounter{example}
    \stepcounter{Mycounter}
    {\noindent \bf \theexample:\ }}

\newcounter{remark}[section]
\renewcommand{\theremark}{{Remark \thesection.\arabic{remark}}}
\newcommand{\remark}{%
    \setcounter{remark}{\value{Mycounter}}
    \refstepcounter{remark}
    \stepcounter{Mycounter}
    {\noindent \bf \theremark:\ }}

\newcounter{problem}[section]
\newcounter{question}[section]
\makeatletter

\@addtoreset{equation}{section} \@addtoreset{footnote}{section} \makeatother

\def\blacksquare{\hbox{\vrule width 5pt height 5pt depth 0pt}}
\def\endproof{\blacksquare}

\begin{document}
\begin{center}
{\LARGE\bf
Existence of HKT metrics on hypercomplex manifolds of
real dimension 8\\[4mm]
}

Gueo Grantcharov, Mehdi Lejmi, Misha Verbitsky\footnote{Gueo 
Grantcharov is supported by a grant from the Simons 
Foundation (\#246184), and Misha Verbitsky is
partially supported by RSCF grant 14-21-00053
within AG Laboratory NRU-HSE.}

\end{center}

{\small \hspace{0.10\linewidth}
\begin{minipage}[t]{0.85\linewidth}
{\bf Abstract} \\
A hypercomplex manifold $M$ is a manifold equipped with three
complex structures satisfying quaternionic relations.
Such a manifold admits a canonical torsion-free connection
preserving the quaternion action, called Obata connection.
A quaternionic Hermitian metric is a Riemannian metric on
which is invariant with respect to unitary quaternions.
Such a metric is called HKT if it is locally obtained
as a Hessian of a function averaged with quaternions.
HKT metric is a natural analogue of a K\"ahler metric on a complex
manifold. We push this analogy further, proving a quaternionic
analogue of Buchdahl-Lamari's theorem for complex surfaces.
Buchdahl and Lamari have shown that a complex surface $M$
admits a Kahler structure iff $b_1(M)$ is even. We show that
a hypercomplex manifold $M$ with Obata holonomy $SL(2,{\mathbb H})$
admits an HKT structure iff $H^{0,1}(M)=H^1({\cal O}_M)$ is even.
\end{minipage}
}

\tableofcontents


\section{Introduction}


\subsection{Hypercomplex manifolds: definition and examples}

Hypercomplex manifolds are the closest quaternionic
counterparts of complex manifolds. They were much studied
by physicists during 1980-ies and 1990-ies, but their
mathematical properties still remain a puzzle.
One obstacle comes from the fact that compact hypercomplex manifold
are non-K\"ahler (unless they are hyperk\"ahler;
see \cite{_Verbitsky:Kahler_HKT_}).
Hypercomplex manifolds appear to be one of
the more-studied and better understood classes of non-K\"ahler
manifolds, which in bigger
generality remain mysterious.
There are many interesting examples of
hypercomplex manifolds and many general theorems,
especially about manifolds admitting HKT-metrics
(\ref{_HKT_Definition_})  or with trivial canonical bundle
(Subsection \ref{_SL(n,H)_Subsection_}).

\hfill

\definition
 Let $M$ be a smooth
manifold equipped with endomorphisms
$I, J, K:\; TM\arrow TM$, satisfying the quaternionic relation
$I^2=J^2=K^2=IJK=-\Id.$  Suppose that $I$, $J$, $K$ are
integrable almost-complex structures. Then $(M, I, J, K)$
is called {\bf a  hypercomplex manifold}.

\hfill

\theorem (Obata, 1955, \cite{_Obata_})\\
On any hypercomplex manifold there exists
a unique torsion-free connection $\nabla$,
called {\bf Obata connection}, such that
$\nabla I = \nabla J =\nabla K$.
\endproof

\hfill

\remark
The holonomy of Obata connection lies in $GL(n, {\Bbb H})$.

\hfill

\remark A torsion-free connection
$\nabla$ on $M$ with $\Hol(\nabla)\subset GL(n, {\Bbb H})$
defines a hypercomplex structure on $M$.

\hfill

\example A {\bf  Hopf surface}
$M={\Bbb H}\backslash 0/ \Z\cong S^1 \times S^3$.
The holonomy of Obata connection $\Hol(M)=\Z$.

\hfill

\example Compact holomorphically symplectic K\"ahler manifolds are
hyperk\"ahler (by Calabi--Yau theorem), hence hypercomplex.
Here $\Hol(M)\subset Sp(n)$ (this holonomy property
is equivalent to being hyperk\"ahler).

\hfill

\proposition
A compact hypercomplex manifold $(M,I,J,K)$
with $(M,I)$ of K\"ahler type
also admits a hyperk\"ahler structure.

{\bf Proof:} \cite[Theorem 1.4]{_Verbitsky:Kahler_HKT_}. \endproof

\hfill

\remark  In quaternionic dimension 1,
compact hypercomplex manifolds are classified
by C. P. Boyer (\cite{_Boyer:note_on_hh_}).
This is the complete list: torus, K3 surface,
Hopf surface.

\hfill

\example
The Lie groups
\begin{align*}
 &  SU(2l+1), \ \ \ \ \ T^1 \times SU(2l), \ \ \  T^l \times SO(2l+1),\\
 & T^{2l}\times SO(4l),  \ \ \ T^l \times Sp(l), \  \ \ \ \ T^2 \times E_6,\\
 & T^7\times E^7, \ \ \ \  \
T^8\times E^8, \ \ \ \ \  \  T^4\times F_4, \ \ \ T^2\times G_2,
\end{align*}
admit a left-invariant hypercomplex structure
(\cite{_SSTvP_}, \cite{_Joyce_}).
Obata holonomy of these manifolds (and other homogeneous
hypercomplex manifolds constructed by Joyce) is unknown,
but most likely it is maximal, that is, equal to $GL(n,\H)$

\hfill

\theorem (Soldatenkov, \cite{_Soldatenkov:SU(2)_}) \\
Holonomy of Obata connection on
$SU(3)$ is $GL(2, {\Bbb H})$. \endproof

\hfill

An important subgroup of $GL(n, {\Bbb H})$ is its commutator $SL(n, {\Bbb H})$. In the standard representation by real matrices this is the subgroup matrices with determinant one. As noted in \cite{_Hitchin_}  it is also isomorphic to one of the real forms of $SL(2n,{\Bbb C})$ denoted by $SU^*(2n)$ in \cite{_Helgason_}. In the present note we focus on manifolds with holonomy in $SL(n, {\Bbb H})$.

\hfill

\example Many {\bf nilmanifolds} (quotients of a nilpotent
Lie group by a cocompact lattice) admit hypercomplex structures.
In this case, $\Hol(M)\subset SL(n, {\Bbb H})$
(\cite{_BDV:nilmanifolds_}).

\subsection{Main result: existence of HKT-metrics on
  $SL(2,{\Bbb H})$-manifolds}

\definition
Let $(M,I,J,K)$ be a hypercomplex manifold, and $g$ a Riemannian
metric. We say that $g$ is {\bf quaternionic Hermitian} if $I,J,K$ are
orthogonal with respect to $g$.

\hfill

\claim Quaternionic Hermitian metrics always exist.

\hfill

{\bf Proof:} Take any Riemannian
metric $g$ and consider its average $\Av_{SU(2)}g$
with respect to $SU(2)\subset {\Bbb H}^*$. \endproof

\hfill

Given a quaternionic Hermitian metric $g$ on $(M,I,J,K)$,
consider its Hermitian forms
\[ \omega_I(\cdot,\cdot)=g(\cdot,I\cdot),\omega_J, \omega_K
\]
(real, but {\em not closed}).
Then $\Omega=\omega_J+\1\omega_K$ is of Hodge
type (2,0) with respect to $I$.

\hfill

\remark If $d\Omega=0,$ the manifold
$(M,I,J,K,g)$ is hyperk\"ahler
(this is one of the definitions of a hyperk\"ahler manifold;
see \cite{_Besse:Einst_Manifo_}).

\hfill

\definition
(Howe, Papadopoulos, \cite{_Howe_Papado_})\\
Let $(M,I,J,K)$ be a hypercomplex manifold,
$g$ a quaternionic Hermitian metric, and
$\Omega=\omega_J+\1\omega_K$
the corresponding $(2,0)$-form. We say that $g$ is HKT
(``hyperk\"ahler with torsion'') if $\6\Omega=0.$

\hfill

\remark
HKT-metrics play in hypercomplex geometry the same role
as K\"ahler metrics play in complex geometry.
\begin{itemize}
\item They admit a smooth potential (locally; see~\cite{_BS:Pot_}).
There is a notion of an ``HKT-class'' (similar to K\"ahler
class) in a certain finite-dimensional
cohomology group, called {\bf Bott--Chern cohomology}
group (Subsection~\ref{_BC_Subsection_}). Two metrics in the same
HKT-class differ by a potential, which is a function.

\item
When $(M,I)$ has trivial canonical bundle,
a version of Hodge theory is established (\cite{_Verbitsky:HKT_}),
giving an ${\goth {sl}}(2)$-action on holomorphic
cohomology $H^*(M, \calo_{(M,I)})$ and analogue of Hodge
decomposition and $dd^c$-lemma.

\item
Originally, it was conjectured that all
hypercomplex manifolds are HKT. The first counterexample
to that assertion is due to Fino and Grantcharov
(\cite{_Fino_Gra_}); for more
examples of non-HKT manifolds, see
\cite{_BDV:nilmanifolds_} and \cite{_Sol_Verb:fibrations_}.
\end{itemize}

The main result of this paper is the following theorem.

\hfill

\theorem\label{_main_intro_Theorem_}
Let $(M,I,J,K)$ be a compact hypercomplex
manifold with Obata holonomy in $SL(2, {\Bbb H})$.
Then  $M$ is HKT if and only if $\dim H^1(\calo_{(M,I)})$ is even.

\hfill

{\bf Proof:} \ref{_HKT_from_even_Corollary_}. \endproof

\hfill

\remark
Using the Hodge decomposition on $H^*(\calo_{(M,I)})$,
one can show that  $h^1(\calo_{(M,I)})$ is even for any
$SL(n, {\Bbb H})$-manifold admitting an HKT-structure
(\cite[Theorem 10.2]{_Verbitsky:HKT_}).

\subsection{Harvey--Lawson duality argument and Lamari's theorem}

The proof of \ref{_main_intro_Theorem_}
is based on the same arguments as used by Lamari (\cite{_Lamari_}) to prove
that any complex surface with even $b_1$ is K\"ahler.
However, in the hypercomplex case this result is
(surprisingly) much easier to prove than in the complex case.

We need the following version of Hahn--Banach theorem:

\hfill

\theorem\label{_H_B_Separation_Theorem_}
(Hahn--Banach separation theorem, \cite{_Schaefer_}) \\
Let $V$ be a locally convex topological
vector space, $A\subset V$ an open convex subset
of $V$, and $W$ a closed subspace of $V$ satisfying $W\cap A=
\emptyset$. Then, there is a continuous linear functional
$\theta$ on $V$, such that
$\theta\restrict A>0$ and $\theta\restrict W = 0$.
\endproof

\hfill

As an illustration, we state the original
Harvey--Lawson duality theorem, which is used as a template
for many other similar arguments, developed since then.

\hfill

\theorem (Harvey, Lawson, \cite{_Harvey_Lawson:Intrinsic_})\\
Let $M$ be a compact complex non-K\"ahler manifold.
Then there exists a positive $(n-1, n-1)$-current
$\xi$ which is a $(n-1,n-1)$-part of an exact current.

\hfill

{\bf Idea of a proof:} Hahn--Banach separation theorem
is applied to the set $A$ of strictly positive $(1,1)$-forms,
and the set $W$ of closed $(1,1)$-forms, obtaining a current
$\xi\in {\mathcal{D}}^{n-1,n-1}(M)= \Lambda^{1,1}(M)^*$ positive on $A$
(that is, positive) and vanishing on $W$. The later condition
(after some simple cohomological manipulations) becomes
``$(n-1,n-1)$-part of an exact current''.
\endproof

\hfill

This approach was further developed some 15 years later by
Buchdahl and Lamari, giving the following theorem.

\hfill

\theorem (Buchdahl--Lamari, \cite{_Buchdahl:surfaces_,_Lamari_})\\
Let $M$ be a compact complex surface. Then $M$ is K\"ahler
if and only if $b_1(M)$ is even.

\hfill

This theorem was known since mid-1980-ies, but its proof
was based on Kodaira classification of complex surfaces,
taking hundreds (if not thousands) of pages and a
complicated result of Siu, who proved that all K3 surfaces
are K\"ahler, and Buchdahl--Lamari (in two independent papers,
\cite{_Buchdahl:surfaces_} and \cite{_Lamari_})
gave a direct proof.

\hfill

{\bf  Scheme of Lamari's proof:}

{\bf  Step 1:} Evenness of $b_1(M)$ is equivalent to
$dd^c$-lemma.

{\bf  Step 2:} Using regularization of positive currents
(\cite{_Demailly_}), one proves that existence of {\bf K\"ahler current}
(positive, closed current $\xi$, such that $\xi -\omega$ is positive
for some Hermitian form $\omega$) is equivalent to existence of a K\"ahler form.

{\bf  Step 3:} Existence of a K\"ahler current is equivalent
to non-existence of a positive current $\xi$ which is a limit of
$dd^c$-closed positive forms and equal to an $(1,1)$-part
of an exact current.

{\bf  Step 4:} Non-existence of such $\xi$ is
implied by $dd^c$-lemma.

\hfill

We are lucky that
for HKT-manifolds the regularization of currents is not necessary and
$dd^c$-lemma
(or, more precisely, its quaternionic analogue)
is the only non-trivial step


\section{Hypercomplex manifolds: basic notions}


\subsection{HKT-manifolds}

The notion of an HKT-manifold
was introduced by the physicists,
but it proved to be immensely useful in mathematics.

A {\bf hypercomplex manifold} is a manifold equipped with
almost-complex structure operators $I, J, K:\; TM \arrow TM$,
integrable and satisfying the standard quaternionic relations
$I^2=J^2=K^2= IJK = -\Id_{TM}$.

This gives a quaternionic algebra action on $TM$;
the group $Sp(1) \cong SU(2)$ of unitary quaternions
acts on all tensor powers of $TM$ by multiplicativity.

A {\bf quaternionic Hermitian structure}
 on a hypercomplex manifold is an $SU(2)$-invariant Riemannian
metric. Such a metric gives a reduction of the structure
group of $M$ to $Sp(n)= U(n, {\Bbb H})$.

With any quaternionic Hermitian structure on $M$
one associates a non-degenerate $(2,0)$-form
$\Omega \in \Lambda^{2,0}_I(M)$, as
follows.\footnote{$\Lambda^*(M)$ denotes the bundle
of differential forms, and $\Lambda^*(M)=\oplus_{p,q}\Lambda^{p,q}_I(M)$
its Hodge decomposition, taken with respect to the
complex structure $I$ on $M$.}
Consider the differential forms
\begin{equation}\label{_three_2_forms_qH_Equation_}
\omega_I(\cdot, \cdot) := g(\cdot, I\cdot), \ \ \omega_J(\cdot, \cdot) := g(\cdot, J\cdot)
, \ \ \omega_K(\cdot, \cdot) := g(\cdot, K\cdot).
\end{equation}
It is easy to check that the form $\Omega:= \omega_J + \1\Omega_K$
is of Hodge type $(2,0)$ with respect to $I$.

If the form $\Omega$ is closed, one has
$d\omega_I = d \omega_J = d\omega_K =0$, and the
manifold $(M,I,J,K, g)$ is called
{\bf hyperk\"ahler} (\cite{_Besse:Einst_Manifo_}).
The hyperk\"ahler condition is very restrictive.

\hfill

\definition\label{_HKT_Definition_}
A hypercomplex, quaternionic Hermitian manifold \\
$(M,I,J,K, g)$ is called {\bf an HKT-manifold}
(hyperk\"ahler with torsion) if $\6\Omega=0$, where
$\6$ denotes the $(1,0)$-part of the differential with respect to $I$.
In other words, a manifold is HKT
if $d\Omega\in \Lambda_I^{2,1}(M)$.

The form $\Omega \in \Lambda_I^{2,0}(M)$
is called {\bf an HKT-form} on $(M,I,J,K)$.

\hfill

\remark\label{_Herm_via_(2,0)_Remark_}
The quaternionic Hermitian form $g$ can be easily
reconstructed from $\Omega$. Indeed, for any
$x, y \in T^{1,0}_I(M)$, one has
\[
2 g(x, \bar y)= \Omega(x, J(\bar y)),
\]
as a trivial calculation implies.

\hfill

Let $(M,I,J,K)$ be a hypercomplex manifold.
We extend \[ J:\; \Lambda^1(M) \arrow \Lambda^1(M)\]
to $\Lambda^*(M)$ by multiplicativity. Recall that
\[ J(\Lambda^{p,q}_I(M))=\Lambda^{q,p}_I(M), \]
because $I$ and $J$ anticommute on $\Lambda^1(M)$.
Denote by
\[ \6_J:\;  \Lambda^{p,q}_I(M)\arrow \Lambda^{p+1,q}_I(M)
\]
the operator $J^{-1}\circ \bar\6 \circ J$, where
$\bar\6:\;  \Lambda^{p,q}_I(M)\arrow \Lambda^{p,q+1}_I(M)$
is the standard Dolbeault operator on $(M,I)$, that is, the
$(0,1)$-part of the de Rham differential.
Since $\bar\6^2=0$, we have $\6_J^2=0$.
In \cite{_Verbitsky:HKT_}, it was shown that $\6$ and $\6_J$
anticommute:
\begin{equation}\label{_commute_6_J_6_Equation_}
\{\6_J, \6 \}=0.
\end{equation}
The pair of anticommuting differentials $\6, \6_J$
is a hypercomplex counterpart to the pair $d, d^c:= I d I^{-1}$
of differentials on a complex manifold.

\subsection{An introduction to $SL(n, {\Bbb H})$-geometry}
\label{_SL(n,H)_Subsection_}

As Obata has shown (\cite{_Obata_}), a hypercomplex
manifold $(M,I,J,K)$ admits a necessarily unique
torsion-free connection, preserving $I,J,K$. The converse
is also true: if a manifold $M$ equipped with an action of
${\Bbb H}$ admits a torsion-free connection preserving the
quaternionic action, it is hypercomplex. This implies that
a hypercomplex structure on a manifold can be defined as a
torsion-free connection with holonomy in $GL(n, {\Bbb
  H})$. This connection is called {\bf the Obata
  connection}.

Connections with restricted holonomy are one of the
central notions in Riemannian geometry, due to Berger's
classification of irreducible holonomy of Riemannian
manifolds. However, a similar classification exists for
general torsion-free connections
(\cite{_Merkulov_Sch:long_}). In the
Merkulov--Schwachh\"ofer list, only three subroups of
$GL(n, {\Bbb H})$ occur. In addition to the compact group
$Sp(n)$ (which defines hyperk\"ahler geometry), also
$GL(n, {\Bbb H})$ and its commutator $SL(n, {\Bbb H})$
appear, corresponding to hypercomplex manifolds and
hypercomplex manifolds with trivial determinant bundle,
respectively. Both of these geometries are interesting,
rich in structure and examples, and deserve detailed
study.

It is easy to see that $(M,I)$ has holomorphically trivial
canonical bundle, for any $SL(n, {\Bbb H})$-manifold
$(M,I, J, K)$ (\cite{_Verbitsky:canoni_}). For a
hypercomplex manifold with trivial canonical bundle
admitting an HKT-metric, a version of Hodge theory was
constructed (\cite{_Verbitsky:HKT_}). Using this result,
it was shown that a compact hypercomplex manifold with
trivial canonical bundle has holonomy in $SL(n,{\Bbb H})$,
if it admits an HKT-structure (\cite{_Verbitsky:canoni_}).

In \cite{_BDV:nilmanifolds_}, it was shown that holonomy
of all hypercomplex nilmanifolds lies in $SL(n, {\Bbb
  H})$. Many  working examples of hypercomplex manifolds
are in fact nilmanifolds, and by this result they all
belong to the class of $SL(n, {\Bbb H})$-manifolds.

The $SL(n, {\Bbb H})$-manifolds were studied in
\cite{_AV:Calabi_} and \cite{_Verbitsky:skoda.tex_}.
On such manifolds the quaternionic Dolbeault
complex is identified with a part of de Rham complex
(\ref{_V_main_Proposition_}), making it possible
to write a quaternionic version of the Monge-Ampere equation
(\cite{_AV:Calabi_}),
and to use quaternionic linear algebra to study positive
currents on hyperk\"ahler manifolds (\cite{_Verbitsky:skoda.tex_}).
 Under this identification, ${\Bbb H}$-positive forms
 become positive in the usual sense, and $\6$,
 $\6_J$-closed or exact forms become $\6, \bar\6$-closed
 or exact (\ref{_V_main_Proposition_}, (iv)). This linear-algebraic
 identification is especially useful in the study of the
 quaternionic Monge-Amp\`ere equation
(\cite{_AV:Calabi_}).

\hfill

One of the main subjects of the present paper
is a quaternionic version of the $dd^c$-lemma,
called ``$\6\6_J$-lemma''.

\hfill

\theorem
Let $M$ be a compact $SL(n,{\Bbb H})$-manifold admitting
an HKT metric, and $\eta$ a $\6_J$-closed, $\6$-exact
$(p,0)$-form. Then $\eta$ lies in the image of
$\6\6_J$.

{\bf Proof:}
In \cite[Theorem 10.2]{_Verbitsky:HKT_}, it was shown that
for any HKT-manifold, the Laplacian $\Delta_\6:= \6\6^*+\6^*\6$ on
$\Lambda^{p,0}(M)\otimes K_M^{1/2}$ can be written as
$\Delta_\6=\{\6,\{\6_J,\Lambda_\Omega\}\}$, where $\{\cdot,\cdot\}$
denotes the anticommutator. Then $\Delta_\6\eta=\6\6_J\Lambda_\Omega\eta$.
However, since $\eta$ is exact, it is orthogonal to the kernel of
$\Delta_\6$, giving $\eta=G\Delta_\6 \eta$, where $G$ is the corresponding
Green operator. This gives
\[
\eta=G\Delta_\6 \eta=G \6\6_J\Lambda_\Omega\eta= \6\6_JG\Lambda_\Omega\eta.
\]
However, on $SL(n,{\Bbb H})$-manifold, the canonical bundle is trivial,
and this result can be applied to any $\eta \in \Lambda^{p,0}(M)$.
\endproof


\section{Quaternionic Dolbeault complex on a hypercomplex manifold}


\subsection{Quaternionic Dolbeault complex: a definition}
\label{_qD_Subsection_}

It is well-known that any irreducible representation
of $SU(2)$ over $\C$ can be obtained as a symmetric power
$\Sym^i(V_1)$, where $V_1$ is a fundamental 2-dimensional
representation. We say that a representation $W$
{\bf has weight $i$} if it is isomorphic to $\Sym^i(V_1)$.
A representation is said to be {\bf pure of weight $i$}
if all its irreducible components have weight $i$.

\hfill

\remark\label{_weight_multi_Remark_}
The Clebsch--Gordan formula (see \cite{_Humphreys_})
claims that the weight is {\em multiplicative},
in the following sense: if $i\leq j$, then
\[
V_i\otimes V_j = \bigoplus_{k=0}^i V_{i+j-2k},
\]
where $V_i=\Sym^i(V_1)$ denotes the irreducible
representation of weight $i$.

\hfill

Let $M$ be a hypercomplex  manifold,
$\dim_{\Bbb H}M=n$.
There is a natural multiplicative action of $SU(2)\subset
{\Bbb H}^*$ on $\Lambda^*(M)$, associated with the
hypercomplex structure.

\hfill

Let $V^i\subset \Lambda^i(M)$ be a maximal
$SU(2)$-invariant subspace of weight $<i$.
The space $V^i$ is well defined, because
it is a sum of all irreducible representations
$W\subset \Lambda^i(M)$ of weight $<i$.
Since the weight is multiplicative
(\ref{_weight_multi_Remark_}), $V^*= \bigoplus_i V^i$
is an ideal in $\Lambda^*(M)$.

It is easy to see that the de Rham differential
$d$ increases the weight by 1 at most. Therefore,
$dV^i\subset V^{i+1}$, and $V^*\subset \Lambda^*(M)$
is a differential ideal in the de Rham DG-algebra
$(\Lambda^*(M), d)$.

\hfill

\definition\label{_qD_Definition_}
Denote by $(\Lambda^*_+(M), d_+)$ the quotient algebra
$\Lambda^*(M)/V^*$.
It is called {\bf the quaternionic Dolbeault algebra of
  $M$}, or {\bf the quaternionic Dolbeault complex}
(qD-algebra or qD-complex for short).

\hfill

\remark
The complex $(\Lambda^*_+(M), d_+)$
was constructed earlier by Capria and Salamon,
(\cite{_Capria-Salamon_}) in a different (and much
more general) situation, and much studied since then.

\subsection{The Hodge decomposition of the quaternionic
  Dolbeault complex}
\label{_Hodge_on_qD_Subsection_}

The Hodge bigrading is compatible with the weight decomposition
of $\Lambda^*(M)$, and gives a Hodge decomposition
of $\Lambda^*_+(M)$ (\cite{_Verbitsky:HKT_}):
\[
\Lambda^i_+(M) = \bigoplus_{p+q=i}\Lambda^{p,q}_{+,I}(M).
\]
The spaces $\Lambda^{p,q}_{+,I}(M)$
are the weight spaces for a particular choice of a Cartan
subalgebra in $\goth{su}(2)$. The $\goth{su}(2)$-action
induces an isomorphism of the weight spaces
within an irreducible representation. This
gives the following result.

\hfill

\proposition \label{_qD_decompo_expli_Proposition_}
Let $(M,I,J,K)$ be a hypercomplex manifold and
\[
\Lambda^i_+(M) = \bigoplus_{p+q=i}\Lambda^{p,q}_{+,I}(M)
\]
the Hodge decomposition of qD-complex defined above.
Then there is a natural isomorphism
\begin{equation}\label{_qD_decompo_Equation_}
\Lambda^{p,q}_{+,I}(M)\cong \Lambda^{p+q,0}(M,I).
\end{equation}

{\bf Proof:} See \cite{_Verbitsky:HKT_}. \endproof

\hfill

This isomorphism is compatible with a natural algebraic
structure on \[ \bigoplus_{p+q=i}\Lambda^{p+q,0}(M,I),\]
and with the Dolbeault differentials, in the following
way.

\hfill

Consider the quaternionic Dolbeault complex
$(\Lambda^*_+(M), d_+)$ constructed in Subsection
\ref{_qD_Subsection_}. Using the Hodge bigrading, we can
decompose this complex, obtaining a bicomplex
\[
\Lambda^{*, *}_{+,I}(M) \xlongrightarrow{d^{1,0}_{+,I}, d^{0,1}_{+,I}}
\Lambda^{*, *}_{+,I}(M)
\]
where $d^{1,0}_{+,I}$,  $d^{0,1}_{+,I}$ are the Hodge components of
the quaternionic Dolbeault differential $d_+$, taken with
respect to $I$.

\hfill

\theorem\label{_bico_ide_Theorem_}
Under the multiplicative isomorphism
\[
\Lambda^{p,q}_{+,I}(M)\cong \Lambda^{p+q,0}(M,I)
\]
constructed in \ref{_qD_decompo_expli_Proposition_},
$d^{1,0}_+$
corresponds to $\6$ and $d^{0,1}_+$
to $\6_J$:
\begin{equation}\label{_bicomple_XY_Equation}
\begin{minipage}[m]{0.85\linewidth}
{\tiny $
\xymatrix @C+1mm @R+10mm@!0  {
  && \Lambda^0_+(M) \ar[dl]^{d^{0,1}_+} \ar[dr]^{d^{1,0}_+}
   && && && \Lambda^{0,0}_I(M) \ar[dl]^{\6} \ar[dr]^{ \6_J}
   &&  \\
 & \Lambda^{1,0}_+(M) \ar[dl]^{d^{0,1}_+} \ar[dr]^{d^{1,0}_+} &
 & \Lambda^{0,1}_+(M) \ar[dl]^{d^{0,1}_+} \ar[dr]^{d^{1,0}_+}&&
\text{\large $\cong$} &
 &\Lambda^{1,0}_I(M)\ar[dl]^{ \6} \ar[dr]^{ \6_J}&  &
 \Lambda^{1,0}_I(M)\ar[dl]^{ \6} \ar[dr]^{ \6_J}&\\
 \Lambda^{2,0}_+(M) && \Lambda^{1,1}_+(M)
   && \Lambda^{0,2}_+(M)& \ \ \ \ \ \ & \Lambda^{2,0}_I(M)& &
\Lambda^{2,0}_I(M) & &\Lambda^{2,0}_I(M) \\
}
$
}
\end{minipage}
\end{equation}
Moreover, under this isomorphism,
$\omega_I\in \Lambda^{1,1}_{+,I}(M)$ corresponds to
$\Omega\in\Lambda^{2,0}_I(M)$.

\hfill

{\bf Proof:} See \cite{_Verbitsky:HKT_} or \cite{_Verbitsky:qD_}.
 \endproof

\subsection{Positive $(2,0)$-forms on hypercomplex
  manifolds}
\label{_posi_2,0-forms_Subsection_}

The notion of positive $(2p,0)$-forms on hypercomplex
manifolds (sometimes called q-positive, or ${\Bbb H}$-positive)
was developed in \cite{_V:reflexive_} and
\cite{_Alesker_Verbitsky_HKT_} (see also
\cite{_AV:Calabi_}
and \cite{_Verbitsky:skoda.tex_}).
For our present purposes, only $(2,0)$-forms are interesting,
but everything can be immediately generalized to a
general situation

\hfill

Let $\eta\in \Lambda^{p,q}_I(M)$ be a
differential form. Since $I$ and $J$ anticommute,
$J(\eta)$ lies in $\Lambda^{q,p}_I(M)$.
Clearly, $J^2\restrict {\Lambda^{p,q}_I(M)}=(-1)^{p+q}$.
For $p+q$ even, $J\restrict {\Lambda^{p,q}_I(M)}$
is an anticomplex involution, that is, a real structure
on $\Lambda^{p,q}_I(M)$.
A form
$\eta \in \Lambda^{2p,0}_I(M)$ is called {\bf real} if
$J(\bar\eta)=\eta$. We denote real forms in $\Lambda^{2p,0}_I(M)$ by
$\Lambda^{2p,0}_{\mathbb{R}}(M,I).$

For a real $(2,0)$-form $\eta$,
\[ 
   \eta\left(x, J(\bar x))\right)=
   \bar \eta\left(J(x), J^2 (\bar x)\right)=
 \bar \eta\left(\bar x, J(x)\right),
\] 
for any $x \in T^{1,0}_I(M)$.
{}From a definition of a real form,
we obtain that the scalar $\eta\left(x, J(\bar x)\right)$
is always real.

\hfill

\definition\label{positive_form}
A real $(2,0)$-form $\eta$ on a hypercomplex manifold
is called {\bf positive} if
$\eta\left(x, J(\bar x)\right)\geq 0$
for any $x \in T^{1,0}_I(M)$, and {\bf strictly positive}
if this inequality is strict, for all $x\neq 0$.

\hfill

An HKT-form $\Omega\in \Lambda^{2,0}_I(M)$ of any
HKT-structure is strictly positive, as follows from
\ref{_Herm_via_(2,0)_Remark_}. Moreover, HKT-structures
on a hypercomplex manifold are in one-to-one
correspondence with closed, strictly positive
$(2,0)$-forms.

The analogy between K\"ahler forms and HKT-forms can be
pushed further: it turns out that any
HKT-form $\Omega\in \Lambda^{2,0}_I(M)$ has a local
potential $\phi\in C^\infty(M)$, in such a way
that $\6\6_J\phi=\Omega$ (\cite{_BS:Pot_}, \cite{_Alesker_Verbitsky_HKT_}).
Here $\6\6_J$ is a composition of $\6$ and $\6_J$
defined on the quaternionic Dolbeault complex as above
(these operators anticommute).

\subsection{The map ${\cal V}_{p,q}:\;
  \Lambda^{p+q,0}_I(M)\arrow\Lambda^{n+p, n+q}_I(M)$\\
on $SL(n, {\Bbb H})$-manifolds}
\label{_V_p,q_Subsection_}

Let $(M,I,J,K)$ be an $SL(n, {\Bbb H})$-manifold, $\dim_\R M =4n$,
and
\[
  {\cal R}_{p,q}:\; \Lambda^{p+q,0}_I(M)\arrow \Lambda^{p,q}_{I,+}(M)
\]
the isomorphism induced by $\goth{su}(2)$-action
as in \ref{_bico_ide_Theorem_}.
Consider the projection
\begin{equation}\label{_proj_to_+_Equation_}
\Lambda^{p,q}_{I}(M)\arrow
\Lambda^{p,q}_{I,+}(M),
\end{equation}
and let $R:\; \Lambda^{p,q}_{I}(M)\arrow\Lambda^{p+q,0}_I(M)$
denote the composition of \eqref{_proj_to_+_Equation_}
and ${\cal R}_{p,q}^{-1}$.

Let $\Phi$ be a nowhere degenerate
holomorphic section of $\Lambda^{2n,0}_I(M)$. Assume that $\Phi$ is
real, that is, $J(\Phi)=\bar\Phi$, and positive.
Existence of such a form is equivalent to
$\Hol(M) \subset SL(n, {\Bbb H})$ (\cite{_Verbitsky:canoni_}).
It is often convenient to define $SL(n, {\Bbb H})$-structure
by fixing the quaternionic action and the holomorphic
form $\Phi$.

\hfill

Define the map
\[ {\cal V}_{p,q}:\;
  \Lambda^{p+q,0}_I(M)\arrow\Lambda^{n+p, n+q}_I(M)
\]
by the relation
\begin{equation}\label{_V_p,q_via_test_form_Equation_}
{\cal V}_{p,q}(\eta) \wedge \alpha = \eta \wedge R(\alpha)\wedge \bar\Phi,
\end{equation}
for any test form $\alpha \in \Lambda^{n-p, n-q}_I(M)$.

\hfill

The map ${\cal V}_{p,p}$ is especially remarkable,
because it maps closed, positive
$(2p,0)$-forms to closed, positive $(n+p, n+p)$-forms,
as the following proposition implies.

\hfill

\proposition\label{_V_main_Proposition_}
Let $(M,I,J,K, \Phi)$ be an $SL(n, {\Bbb H})$-manifold, and
\[ {\cal V}_{p,q}:\;
  \Lambda^{p+q,0}_I(M)\arrow\Lambda^{4n-p, 4n-q}_I(M)
\]
 the map defined above.
Then
\begin{description}
\item[(i)] ${\cal V}_{p,q}(\eta)=
{\cal R}_{p,q}(\eta) \wedge {\cal V}_{0,0}(1)$.
\item[(ii)]  The map ${\cal V}_{p,q}$ is injective, for
  all $p$, $q$.
\item[(iii)] $(\1)^{(n-p)^2}{\cal V}_{p,p}(\eta)$ is real if and
  only $\eta\in\Lambda^{2p,0}_I(M)$ is real,
and weakly positive if and only if $\eta$ is weakly positive.
\item[(iv)] ${\cal V}_{p,q}(\6\eta)= \6{\cal V}_{p-1,q}(\eta)$,
and ${\cal V}_{p,q}(\6_J\eta)= \bar\6{\cal  V}_{p,q-1}(\eta)$.
\item[(v)] ${\cal V}_{0,0}(1) = \lambda {\cal
  R}_{n,n}(\Phi)$, where $\lambda$ is a positive rational number,
depending only on the dimension $n$.
\end{description}

{\bf Proof:} See \cite{_Verbitsky:skoda.tex_},
Proposition 4.2, or \cite{_AV:Calabi_}, Theorem 3.6. \endproof


\section{Quaternionic Gauduchon metrics}


\subsection{Gauduchon metrics}

\definition
A Hermitian metric $\omega$ on a complex $n$ manifold is
called {\bf Gauduchon} if $\6\bar\6\omega^{n-1}=0$.

\hfill

\theorem
Every Hermitian metric on a compact complex manifold
is conformally equivalent to a Gauduchon metric,
which is unique in its conformal class, up to a
constant multiplier.

{\bf Proof:} \cite{_Gauduchon_1984_}. \endproof

\hfill

Gauduchon metrics is one of the very few instruments available
for the study of general non-K\"ahler manifolds, and probably the most
important one.

\subsection{Gauduchon metrics and hypercomplex structures}

Let $g$ be a quaternionic Hermitian metric on
a hypercomplex manifold $M$.
Consider the corresponding $(2,0)$-form
$\Omega:=\omega_J+\1 \omega_K$ defined as in
\ref{_Herm_via_(2,0)_Remark_}. From the definition
of positive (2,0)-forms it follows that
this correspondence is bijective: quaternionic
Hermitian metrics are in (1,1)-correspondence
with positive (2,0)-forms.

\hfill

\definition
A quaternionic Hermitian form $g$ on
a hypercomplex manifold $M$, $\dim_{\Bbb H}M=n$,
is called {\bf quaternionic Gauduchon} if
$\6\6_J\Omega^{n-1}=0$, where $\Omega=\omega_J+\1\omega_K$
is the corresponding positive (2,0)-form.

\hfill

\proposition
Let $(M,I,J,K,\Phi)$ be an $SL(n,{\Bbb H})$-manifold
equipped with a quaternionic Hermitian form $g$, and
\[
|\Phi|^2:=\frac{\Phi\wedge\bar\Phi}{(2^{2n}2n!)^{-1}\omega_I^{2n}}
\]
Then the following conditions are equivalent.
\begin{description}
\item[(i)] $g$ is quaternionic Gauduchon.
\item[(ii)] The Hermitian metric $|\Phi|^{-1}g$ is Gauduchon on $(M,I)$.
\item[(iii)] The Hermitian metric $|\Phi|^{-1}g$ is Gauduchon with respect
to any of the induced complex structures $L=aI+bJ+cK$
\end{description}


{\bf Proof:} The equivalence (i) $\Leftrightarrow$
(ii) follows from
\begin{equation*}
{\cal V}_{n-1,n-1}(\Omega^{n-1})=|\Phi|^{-1}\omega_I^{2n-1},
\end{equation*}
proven in \cite{_Gra_Verb_} (the formula in the proof
of Theorem 6.4).
So, using \ref{_V_main_Proposition_} (iv), we have that $ {\cal V}_{n,n}(\6\6_J\Omega^{n-1})=\6\bar\6 (|\Phi|^{-1}\omega_I^{2n-1})$.

\endproof

\hfill

\corollary
For any $SL(n,{\Bbb H})$-manifold
equipped with a quaternionic Hermitian form, there
exists a unique (up to a constant multiplier) positive
function $\mu$ such that $\mu g$ is quaternionic Gauduchon.
\endproof

\subsection{Surjectivity of $f\arrow \Omega^{n-1}\wedge\6\6_Jf$/}

We are interested in quaternionic  Gauduchon forms
because of the following theorem.

\hfill

\theorem\label{_elli_66_J_bije_Theorem_}
Let $(M,I,J,K,\Omega,\Phi)$
be a compact quaternionic Hermitian $SL(n,{\Bbb H})$-manifold. Assume that
$\Omega$ is quaternionic Gauduchon. Consider the map
$D:\; C^\infty (M) \arrow \Lambda^{4n}(M)$,
\[ D(f)= \6\6_Jf \wedge \Omega^{n-1}\wedge \Phi.\]
Then $D$ induces a bijection between $C^\infty (M)/\const$
and the space of exact $4n$-forms on $M$.

\hfill

{\bf  Proof: Step 1:} Clearly, $D$ is elliptic, and
has index 0, because it has the same symbol as Laplacian,
which is self-adjoint.

{\bf  Step 2:} E. Hopf maximum principle
(\cite{_Gilbarg_Trudinger_}) implies that $\ker D=\const$.
Therefore, $\coker D$ is 1-dimensional. It remains to show that
$\im D$ consists of exact $4n$-forms.

{\bf  Step 3:}
\[
 \int_M \6\6_Jf \wedge \Omega^{n-1}\wedge \bar\Phi=
 - \int_M f \wedge \6\6_J(\Omega^{n-1})\wedge \bar\Phi=0
\]
because $\Omega$ is quaternionic Gauduchon. This implies that
all forms in $\im D$ are exact. Converse is also true,
because $\codim \im D=1$.
\endproof


\section{Quaternionic Aeppli and Bott--Chern cohomology}


Throughout this section, $(M,I,J,K,g)$ is a compact hypercomplex manifold equipped with a quaternionic Hermitian metric $g$. Recall that $\{\6,\6_J\}=0$.

\subsection{Quaternionic Bott--Chern cohomology}
\label{_BC_Subsection_}

Define $H^{p,0}_{BC}(M)$ to be the
group $$H^{p,0}_{BC}(M)= \frac{\{ \phi\in
  {\Lambda}^{p,0}_I(M)\, |\, \6 \phi =\6_J\phi=
  0\}}{\6\6_J{\Lambda}^{p-2,0}_I(M)}.$$

\theorem
The group $H^{p,0}_{BC}(M)$ is finite dimensional.

\hfill

{\bf Proof:}
We consider the following operator $${\Delta}_{BC}=\6^\ast\6+\6_J^\ast\6_J+\6\6_J\6_J^\ast\6^\ast+\6_J^\ast\6^\ast\6\6_J+\6_J^\ast\6\6^\ast\6_J+\6^\ast\6_J\6_J^\ast\6,$$ acting on
${\Lambda}^{p,0}_I(M)$. Here, $\6^\ast$ (resp. $\6_J^\ast$) is the adjoint of $\6$ (resp. $\6_J$) with respect to $g.$

We claim that ${\Delta}_{BC}$ is a fourth order self-adjoint elliptic operator.
Using the elliptic theory, we obtain the following decomposition
\begin{eqnarray*}
 {\Lambda}^{p,0}_I(M)&=&\mathcal{H}_{{\Delta}_{BC}}\oplus \im\,{\Delta}_{BC},\\
 &=&\mathcal{H}_{{\Delta}_{BC}}\oplus \im\,\6\6_J\oplus\,( \im\,\6^\ast+\, \im\,\6_J^\ast),
\end{eqnarray*}
where $\mathcal{H}_{{\Delta}_{BC}}=\{ \phi\in{\Lambda}^{p,0}_I(M)\,|\, \6\phi=\6_J\phi=\6_J^\ast\6^\ast\phi=0\}$ is the kernel of ${\Delta}_{BC}.$


Furthermore, for $\phi\in{\Lambda}^{p,0}_I(M)$, we write $\phi=\phi_H+\6\6_J\rho+\6^\ast\alpha+\6_J^\ast\beta$, where $\phi_H\in\mathcal{H}_{{\Delta}_{BC}}.$ Then,
$\6\phi=\6_J\phi=0$ is equivalent to $\6^\ast\alpha+\6_J^\ast\beta=0.$
Thus, we deduce $$\ker \6|_{ {\Lambda}_I^{p,0}(M)}\,\cap\,\ker\6_J|_{ {\Lambda}_I^{p,0}(M)}=\mathcal{H}_{{\Delta}_{BC}}\oplus \im\,\6\6_J.$$\endproof

\subsection{Quaternionic Aeppli cohomology}

In a similar way, we define $H^{p,0}_{AE}(M)$
to be the group $$H^{p,0}_{AE}= \frac{\{ \phi\in
  {\Lambda}_I^{p,0}(M)\, |\, \6\6_J\phi= 0\}}{\6{\Lambda}^{p-1,0}_I(M)+\6_J{\Lambda}^{p-1,0}_I(M)}.$$

\theorem
The group $H^{p,0}_{AE}(M)$ is finite dimensional.

\hfill

{\bf Proof:}
 Here, we consider the operator $${\Delta}_{AE}=\6\6^\ast+\6_J\6_J^\ast+\6\6_J\6_J^\ast\6^\ast+\6_J^\ast\6^\ast\6\6_J+\6\6_J^\ast\6_J\6^\ast+\6_J\6^\ast\6\6_J^\ast,$$ acting on
${\Lambda}_I^{p,0}(M)$.
The operator ${\Delta}_{AE}$ is a fourth order self-adjoint elliptic operator having the same symbol as ${\Delta}_{BC}.$
We have then
\begin{eqnarray*}
 {\Lambda}_I^{p,0}(M)&=&{{\mathcal{H}}}_{{\Delta}_{AE}}\oplus \im\,{\Delta_{AE}},\\
 &=&{{\mathcal{H}}}_{{\Delta}_{AE}}\oplus \im\,\6_J^\ast\6^\ast\oplus\,( \im\,\6+\, \im\,\6_J),
\end{eqnarray*}
where $\mathcal{H}_{{\Delta}_{AE}}=\{ \phi\in{\Lambda}^{p,0}_I(M)\,|\,\6^\ast\phi=\6_J^\ast\phi=\6\6_J\phi=0\}$ is the kernel of ${\Delta}_{AE}.$

%

Moreover, if $\phi\in{\Lambda}^{p,0}_I(M)$ is decomposed as $\phi=\phi_H+\6_J^\ast\6^\ast\rho+\6\alpha+\6_J\beta,$ where $\phi_H\in{{\mathcal{H}}}_{{\Delta}_{AE}},$
then $\6\6_J\phi=0$ is equivalent to $\6_J^\ast\6^\ast\rho=0.$
We obtain that $$\ker \6\6_J|_{ {\Lambda}_I^{p,0}(M)}={{\mathcal{H}}}_{{\Delta}_{AE}}\oplus( \im\,\6+\, \im\,\6_J).$$
\endproof

\hfill

\remark\label{_duality_} The groups $H^{p,0}_{BC}(M)$ and $H^{2n-p,0}_{AE}(M)$ are dual when $M$ is a compact $SL(n,\mathbb{H})$-manifold. Indeed, let $\Phi$ be a nowhere degenerate holomorphic section of $\Lambda^{2n,0}_I(M)$.
We assume also that $\Phi$ is real and positive.
We consider the pairing on $H^{p,0}_{BC}(M)\times H^{2n-p}_{AE}(M)$ given by $$([\alpha],[\beta])\mapsto\int_M\alpha\wedge\beta\wedge\bar\Phi.$$ One can check that this pairing is well defined (recall that $\6\bar\Phi=\6_J\bar\Phi=0$)
and non-degenerate.


\section{$\6\6_J$-lemma in $\dim_{\Bbb H}=2$.}


\definition
Let $(M,I,J,K,\Omega,\Phi)$
be a compact quaternionic Gauduchon $SL(n,{\Bbb H})$-manifold,
and $H^{*,0}_{AE}(M)$, $H^{*,0}_{BC}(M)$ the quaternionic Aeppli
and Bott--Chern cohomology.
Consider the map $\deg : H^{1,0}_{AE}(M)\arrow \C$ putting
$\alpha$ to $\int\6\alpha\wedge \Omega^{n-1}\wedge \bar\Phi$.
Since $\Omega$ is quaternionic Gauduchon, $\deg\alpha$ is independent
from the choice of $\alpha$ in its cohomology class.
We call $\deg$ {\bf the degree map}.

\hfill

\remark
Consider the natural map
$H^{1,0}_{AE}(M)\stackrel{\6}\arrow H^{2,0}_{BC}(M)$.
The kernel of this map consists of all
cohomology classes $\alpha$ such that $\6\alpha=\6\6_J\beta$,
hence the form $\alpha-\6_J\beta$ cohomological to $\alpha$
in $H^{1,0}_{AE}(M)$ is $\6$-closed. We obtain that
the kernel of $H^{1,0}_{AE}(M)\stackrel{\6}\arrow H^{2,0}_{BC}(M)$
is identified with the space
$H^{1,0}_\6(M)=\frac{\ker \6\restrict{\Lambda^{1,0}(M)}}{\im\6}$.

\hfill

\lemma\label{_66_J_for_H^1_Lemma_}
($\6\6_J$-lemma for $H^{1,0}(M)$)\\
Let $(M,I,J,K,\Omega,\Phi)$
be a compact quaternionic Gauduchon $SL(2,{\Bbb H})$-\-ma\-nifold. Let $\theta \in \Lambda^{1,0}_I(M)$ be a $\6_J$-exact,
$\6$-closed form. Then $\theta=0$.

\hfill

{\bf Proof:} Let $\theta= 6_J(f)$. Then $\6\6_J(f)=\6\theta=0$.
However, the map $f\arrow \frac{\6\6_J(F)\wedge \Omega}{\Omega^2}$
is an elliptic operator with vanishing constant term, hence
any function in its kernel is constant by Hopf maximum principle
(\cite{_Gilbarg_Trudinger_}).
\endproof

\hfill

\corollary\label{_injectivity_}
On a compact $SL(2,{\Bbb H})$-manifold, the natural map
\[ H^{1,0}_\6(M)\longrightarrow H^{1,0}_{AE}(M)
\] is injective.

\hfill

\theorem\label{_degree_H^1_AE_Theorem_}
Let $(M,I,J,K,\Omega,\Phi)$
be a compact quaternionic Gauduchon $SL(2,{\Bbb H})$-manifold.
Then, the sequence
\begin{equation}\label{_AE_coho_degree_sequece_exact_Equation_}
0\arrow H^{1,0}_\6(M)\arrow H^{1,0}_{AE}(M)\stackrel{\deg}\arrow \C
\end{equation}
is exact. Moreover, the space $\ker \deg$ is equal to the kernel
of the natural map $H^{1,0}_{AE}(M)\stackrel{\6}\arrow H^{2,0}_{BC}(M)$.

\hfill

{\bf Proof: Step 0:} By \ref{_66_J_for_H^1_Lemma_},
the sequence \eqref{_AE_coho_degree_sequece_exact_Equation_}
is exact in the first term. It remains to prove that
\eqref{_AE_coho_degree_sequece_exact_Equation_} is exact in the
second term and to show that $\ker \deg\ker \6 \restrict{H^{1,0}_{AE}(M)}$.

\hfill

{\bf Step 1:}
Let $\alpha\in \ker \deg$. By
\ref{_elli_66_J_bije_Theorem_}, there exists $f\in C^\infty (M)$ such that
$(\6\alpha + \6\6_J f)\wedge \Omega\wedge \bar\Phi=0$,
equivalently $(\6\alpha + \6\6_J f)\wedge \Omega=0$.
Replacing $\alpha$ by $\alpha+\6_J f$ in the same cohomology
class,  we may assume that $\6\alpha\wedge\Omega=0$.

\hfill

{\bf Step 2:} Since $\6\alpha$ is primitive,
one has $\int_M \6\alpha\wedge \6_J\alpha\wedge \bar\Phi=-\|\6\alpha\|^2$
by a quaternionic version of Hodge--Riemann relations
(\cite[Theorem 6.3]{_Verbitsky:balanced_}).

\hfill

{\bf Step 3:} However,
\[ \|\6\alpha\|^2=
\int_M \6\alpha\wedge \6_J\alpha\wedge \bar\Phi=-
\int_M \6\6_J\alpha\wedge \alpha\wedge \bar\Phi=0,
\] hence
$\6\alpha=0$. This implies that $\ker\deg=\ker \6$.
\endproof

\hfill

The $\6\6_J$-lemma for even $h^1(\calo_M)$ follows directly
from the above theorem.

\hfill

\theorem\label{_6_6_J_lemma_even_coho_Theorem_}
Let $(M,I,J,K,\Phi)$
be a compact $SL(2,{\Bbb H})$-manifold.
Then $\6\6_J$-lemma holds on $\Lambda^{2,0}(M)$
if and only if $h^1(\calo_M)$ is even.

\hfill

{\bf Proof: Step 1:} Clearly, $\6\6_J$-lemma
is equivalent to vanishing of $\6:\; H^{1,0}_{AE}(M)\arrow H^{2,0}_{BC}(M)$,
but the kernel of this map is $H^{1,0}_\6(M)=\ker \deg$
(\ref{_degree_H^1_AE_Theorem_}), hence it suffices
to show that the degree map vanishes
iff $h^{1,0}(\calo_M)$ is even.

{\bf Step 2:} Since $J$ defines the quaternionic
structure on $H^{1,0}_{AE}(M)$, this space is even-dimensional.
Now, from the exact sequence
\[ 0\arrow H^{1,0}_\6(M)\arrow H^{1,0}_{AE}(M)\stackrel{\deg}\arrow \C,\]
we obtain that $\deg=0$ whenever $H^{1,0}_\6(M)$ is even-dimensional.
The space $H^{1,0}_\6(M)$ is complex conjugate to $H^1(\calo_M)$,
hence has the same dimension.
\endproof


\section{Currents in HKT-geometry}


\subsection{Cohomology of currents}

\definition
Let $(M,I,J,K)$ be a hypercomplex manifold. Denote by ${\mathcal D}_{p,q}(M)$ the topological dual to the Fr\'echet space $\Lambda^{p,q}_{I}(M)$. An element $T\in{\mathcal D}_{p,q}(M)$
is called a {\bf current} of bidimension $(p,q)$ and it has a compact support on $M$. Denote by ${\mathcal D}^{p,q}(M)={\mathcal D}_{2n-p,2n-q}(M)$, where $\dim_{\mathbb{H}}M=n.$

\hfill

The complex structure $J$ acts naturally on ${\mathcal D}^{p,q}(M)$ as a map $$J: {\mathcal D}^{p,q}(M)\rightarrow{\mathcal D}^{q,p}(M)$$ in the following way $$(JT)(\phi)=T(J\phi),$$
for $\phi \in \Lambda_I^{2n-q,2n-p}(M)$ with compact support.
The operators $d, \6, \overline{\6}$ are extended in the standard way
using the Stokes theorem, for example
$\6:{\mathcal D}^{p,q}(M)\rightarrow {\mathcal D}^{p+1,q}(M)$
is expressed as $\6T(\phi) = (-1)^{\dim \phi}T(\6\phi),$ where
$\phi \in \Lambda_I^{2n-p-1,2n-q}(M)$. Similarly, we can define $\6_J = J^{-1}\circ\bar\6\circ J$ on ${\mathcal D}^{p,q}(M)$.

\hfill

\definition
A current $T\in{\mathcal D}^{2p,0}(M)$ is called {\bf real} if $J\bar{T}=T$ and we denote real currents by ${\mathcal D}^{2p,0}_{\mathbb{R}}(M)$.







\hfill

The following result is a currents version of
local $\6\6_J$-lemma, due to Banos and Swann in the smooth case.

\hfill

\proposition
Let $T \in{\mathcal D}_{\mathbb{R}}^{2,0}(M)$ be a real $\6$-closed current. Then, locally $T$ can be written in the form $T=\6\6_J\varphi,$
for some real generalized function $\varphi.$

\hfill

{\bf Proof:}
We use essentially the same arguments as in the proof of the main theorem in \cite{_BS:Pot_}. Let $T$ as above. We write $T=T_J+\sqrt{-1}\,T_K.$ Since $T$ is real and $\6$-closed, a straightforward verification shows that $T_I(\cdot,\cdot)=T_J(\cdot,K\cdot)$
is $I$-invariant and that $IdT_I=JdT_J=KdT_K$.

Let $Z= M\times S^2$ be the twistor space of $M$ and we consider the current $\eta\in{\mathcal D}^{0,2}(Z)$ given by $\eta=(T_I)^{(0,2)}_\mathcal{I}$ i.e. the $(0,2)$-part
of $T_I$ with respect to the complex structure $\mathcal{I}|_{(p,\overrightarrow{a})}=aI+bJ+cK$ ($\overrightarrow{a}=(a,b,c)\in\mathbb{R}^3,\,a^2+b^2+c^2=1$). A direct computation shows that $\bar\6_{\mathcal{I}}\eta=0$. By a 1-pseudo-convexity argument and the $\bar\6$-Poincar\'e Lemma (for currents), locally $\eta=\bar\6_{\mathcal{I}}(\alpha+\sqrt{-1}\,\mathcal{I}\,\alpha)$ where $\alpha$ is a real current defined locally in $M$. Hence, the real part of $\eta$ is given by $\frac{1}{2}(d\alpha-\mathcal{I}d\alpha)$. It follows that $d\alpha$ is a closed $I$-invariant current.
Hence, by the $\6\bar\6$-Poincar\'e Lemma (for currents), $T_I=\frac{1}{2}(dd_I\varphi+d_Jd_K\varphi)$ for some real generalized function $\varphi$. By \cite{_Alesker_Verbitsky_HKT_}, this implies that locally $T=\6\6_J\varphi.$
\endproof

\hfill

Using~\ref{positive_form}, we give the following:

\hfill

\definition On a $SL(n,\mathbb{H})$-manifold $(M,I,J,K,\Phi)$,
a current $T\in{\mathcal D}^{2n-2,0}_{\mathbb{R}}(M)$
is said to be {\bf positive} if (locally) $T\wedge\alpha\wedge\bar\Phi$ is a positive measure for any choice of (local) real strictly positive
$(2,0)$-form $\alpha$.

\hfill

\definition
A generalized function is called {\bf plurisubharmonic}
if $\6\6_J\phi$ is a positive (2,0)-current.

\hfill

\theorem\label{_psh_cur_compa_constant_Theorem_}
A plurisubharmonic generalized function
is subharmonic with respect to any
quaternionic Hermitian metric (hence, constant
on any compact hypercomplex manifold).

\hfill

{\bf Proof:}
\cite[Lemma 3.6]{_Harvey_Lawson:Potential_Calibrated_}.
\endproof

\hfill

\remark\label{_currents_duality_}
We consider the group $$H'^{2,0}_{\mathbb{R}}(M)= \frac{\{u \in {\mathcal D}_{\mathbb{R}}^{2,0}(M)\,|\,\6 u=0\}}{\6\6_J{\mathcal D}_{\mathbb{R}}^{0,0}(M)}.$$
Denote by $\mathcal{H}$ the sheaf of real generalized functions
satisfying $\6\6_Jf=0$. By the proof of
\ref{_66_J_for_H^1_Lemma_}, elements of $\mathcal{H}$
satisfy an elliptic equation. Elliptic regularity implies
that all functions in $\mathcal{H}$ are smooth.

The sheaf $\mathcal{H}$ admits two resolutions starting by
\begin{diagram}
0 &\rTo & {\mathcal H}&\rTo & \Lambda_{\mathbb{R}}^{0,0}(M,I)&\rTo^{\6\6_J} &\Lambda^{2,0}_{\mathbb{R}}(M,I)&\rTo^{\6} &\Lambda_{}^{3,0}(M,I)\\
&             &\dTo_{id}         &             &\dTo_{i}                   &                       &\dTo_{i}                                     &                                                       &\dTo_{i}\\
0 &\rTo & {\mathcal H}&\rTo& {\mathcal D}_{\mathbb{R}}^{0,0}(M) &\rTo^{\6\6_J} &{\mathcal D}_{\mathbb{R}}^{2,0}(M)&\rTo^{\6} &{\mathcal D}_{}^{3,0}(M),
\end{diagram}
where $i$ is the inclusion of forms in the space of currents. We deduce that $$H'^{2,0}_{\mathbb{R}}(M)\simeq\frac{\{ \phi\in \Lambda^{2,0}_{\mathbb{R}}(M,I)\,|\,\6 \phi= 0\}}{\6\6_J{\Lambda}_{\mathbb{R}}^{0,0}(M,I)}.$$



Let $$H'^{2n-2,0}_{\mathbb{R}}(M)= \frac{\{u \in {\mathcal D}_{\mathbb{R}}^{2n-2,0}(M)\,|\,\6\6_J u=0\}}{\{\6\eta+\6_JJ^{-1}\bar\eta,\,\, \eta\in{\mathcal D}_{\mathbb{}}^{2n-3,0}(M)\}}.$$ Then, by the same argument in \ref{_duality_}, we deduce that
$H'^{2n-2,0}_{\mathbb{R}}(M)$ and $H'^{2,0}_{\mathbb{R}}(M)$ are dual when $M$ is a compact $SL(n,\mathbb{H})$-manifold.

%


%

\hfill

\subsection{Harvey--Lawson's theorem in HKT-geometry}


Using the Hahn--Banach Separation Theorem
(\ref{_H_B_Separation_Theorem_}), we obtain the following.

\hfill

\theorem\label{_HL_for_HKT_Theorem_}
Let $(M,I,J,K,\Phi)$ be an $SL(n,{\Bbb H})$-manifold.
Then $M$ admits no HKT-metric if and only if
it admits a $\6$-exact, real, positive $(2n-2,0)$-current.

\hfill

{\bf  Proof: Step 1:} Apply Hahn--Banach separation theorem
to the space $A$ of positive, real $(2,0)$-forms and
$W$ of $\6$-closed real (2,0)-forms to obtain a current
$\xi\in \Lambda^{2,0}_\R(M,I)^*$ which is positive on $A$
(hence, real and positive) and vanishes on $W$.
Such a current exists iff $A\cap W=\emptyset$, or, equivalently,
when $M$ is not HKT.

{\bf  Step 2:} Consider the pairing
$\langle \eta, \nu\rangle = \int_M \eta\wedge \nu\wedge
\bar \Phi$ on $(p,0)$-forms. This pairing is compatible
with $\6$ and $\6_J$ and allows one to identify
the currents $\Lambda^{p,0}_\R(M,I)^*$ with
$\Lambda^{n-p,0}_\R(M,I)\otimes C^\infty(M)^*$,
where $C^\infty(M)^*$ denotes generalized functions.
This identification is compatible with $\6$ and $\6_J$,
and cohomology of currents are the same as cohomology of
forms (\ref{_currents_duality_}).

{\bf Step 3:}
Since $\langle \xi, W\rangle =0$, for each $\eta$
one has $0=\langle \xi, \6\eta\rangle =\langle \6\xi,
\eta\rangle$, giving $\6\xi=0$. It remains to show that
the cohomology class of $\xi$ in $H^2_{\6}(\Lambda^{*,0}_I(M))$ vanishes.

{\bf Step 4:} The Serre's duality gives a
non-degenerate pairing $\langle [\xi], [\nu]\rangle\arrow \R$
on cohomology classes in  $H^2_{\6}(\Lambda^{*,0}_I(M))$:
\[
\xi, \nu \arrow \int_M \xi\wedge\nu\wedge\bar\Phi.
\]
 Since
$\langle [\xi], [\nu]\rangle=0$ for each $\6$-closed $\nu$,
the cohomology class of $\xi$ also vanishes.
\endproof

\hfill

\corollary\label{_HKT_from_even_Corollary_}
Let $M$ be a compact $SL(2,{\Bbb H})$-manifold. Then $M$ admits HKT-metric
if and only if $H^1(\calo_{(M,I)})$ is even-dimensional.

\hfill

{\bf Proof:}
Even-dimensionality of $H^1(\calo_{(M,I)})$ for HKT-manifolds
with holonomy in $SL(n,\H)$ follows from
\cite[Theorem 10.2]{_Verbitsky:HKT_}. Conversely,
suppose that $H^1(\calo_{(M,I)})$ is even-dimensional,
but $M$ is not HKT. Then \ref{_HL_for_HKT_Theorem_}
implies that there exists a real, positive, exact (2,0)-current $\xi$.
However, $\xi$ is $\6\6_J$-exact by \ref{_6_6_J_lemma_even_coho_Theorem_},
hence $\xi=\6\6_Jf$, for some $f\in C^\infty (M)$. Such $f$
is a quaternionic plurisubharmonic function, which has to vanish
by \ref{_psh_cur_compa_constant_Theorem_}.
\endproof

\hfill


\section{Examples}


The known examples of manifolds with holonomy $SL(n,{\Bbb
  H})$ are either nilmanifolds
(\cite{_BDV:nilmanifolds_}) or
obtained via the twist construction of A.~Swann \cite{_Swann1_},
which is based on previous examples by D.~Joyce. The later
construction provides also simply-connected examples.  We describe briefly a simplified version of it.

Let $(X, I,J,K,g)$ be a compact hyperk\"ahler
manifold. By definition, an anti-self-dual 2-form on it is
a form which is of type (1,1) with respect to $I$ and $J$
and hence with respect to all complex structures of the
hypercomplex family. Let $\alpha_1,\cdots,\alpha_{4k}$ be
closed 2-forms representing integral
cohomology classes on $X$. Consider the principal $T^{4k}$-bundle $\pi:M\rightarrow X$ with
characteristic classes determined by $\alpha_1,\cdots,\alpha_{4k}$. It admits a connection $A$ given by
$4k$ 1-forms $\theta_i$ such that $d\theta_i=\pi^*(\alpha_i)$. Define an almost-hypercomplex structure on $M$
in the following way: on the horizontal spaces of $A$ we have the pull-backs of $I,J,K$ and on the vertical
spaces we fix a linear hypercomplex structure of the
$4k$-torus. The structures $\cal{I}, \cal{J}, \cal{K}$
on $M$ are extended to act on the cotangent bundle $T^*M$
using the following relations:
\begin{align*}
{\cal I}(\theta_{4i+1})=\theta_{4i+2}, \,\,&{\cal
  I}(\theta_{4i+3})=\theta_{4i+4}, &{\cal
  J}(\theta_{4i+1})=\theta_{4i+3},\,\, &{\cal
  J}(\theta_{4i+2})=-\theta_{4i+4},\\
 {\cal I}(\pi^*\alpha)=\pi^*(I\alpha),\,\, &{\cal
   J}(\pi^*\alpha)=\pi^*(J\alpha),& &
\end{align*}
for any 1-form $\alpha$ on $X$ and $i=0,\cdots,k-1$.

It follows from \cite{_Swann1_} or by direct and easy calculations, that $\cal{I}$ is integrable iff $\alpha_{4i+1}+i \alpha_{4i+2}$ and $\alpha_{4i+3}+ i\alpha_{4i+4}$ are of type $(2,0)+(1,1)$ with respect to $I$ for every $i = 1,\cdots,k$. Similarly $\cal{J}$ is integrable iff $\alpha_{4i+1}+i \alpha_{4i+3}$ and $\alpha_{4i+2}- i\alpha_{4i+4}$ are of type $(2,0)+(1,1)$ with respect to $J$ for every $i = 1,\cdots,k$

Similarly, one can define a
quaternionic Hermitian metric on $M$ from $g$ and a fixed hyperk\"ahler metric on $T^{4k}$ using the splitting
of $T(M)$ in horizontal and vertical subspaces. As A.~Swann \cite{_Swann1_} has shown the structure  has
a holonomy in $SL(n,{\Bbb H})$ and is HKT when all forms $\alpha_l$ are self-dual (of type $(1,1)$ with respect to all structures).

 As a particular case, assume $X$ to be a $K3$ such that there are 3 closed integral forms which define a hyperk\"ahler structure
 and a self-dual integral class, so defining a
 principal $T^4$-bundle $M$ over $X=K3$ with finite
 fundamental group. After passing to a finite cover, we may
 assume that $M$ is simply-connected. These forms satisfy
 the integrability condition above. If $\alpha_2+i
 \alpha_3$ is a (2,0)-form for $I$, then $\pi^*(\alpha_2+i
 \alpha_3) = d(\theta_2+i\theta_3)$ is an exact
 $(2,0)$-form, which defines a positive current in the
 definition of the previous section. Then $M$ can not
 admit any HKT-metric - a fact proven by A. Swann using
 different arguments.



 We can also calculate $\dim(H^1(\calo_{(M,\cal I)})) =
 h^{0,1}_{\cal I}(M)$ and apply~\ref{_main_intro_Theorem_} to decide the
 existence of HKT-structure. One can use the Borel method
 of doubly graded spectral sequence from
 \cite{_Hirzebruch_}, Appendix B, to determine $h^{p,q}$,
 but in our case, its simpler to use a more direct
 approach. The vector fields $X_1, X_2, X_3, X_4$ on $M$
 generated by torus action which are also dual to
 $\theta_i$ are hyperholomorphic, so ${\cal
   L}_{X_i}\circ{\cal I} = {\cal I}\circ {\cal
   L}_{X_i}$. We can also choose  a bundle metric, which
 for the vertical vectors is the flat hyperk\"ahler
 4-torus metric and on the horizontal vectors is a
 pull-back from the hyperk\"ahler metric from the base
 $X=K3$. The horizontal and vertical vectors are
 perpendicular. Such metric is hypercomplex and $X_i$ are
 Killing fields. So, since they also fix the orientation,
 then ${\cal L}_{X_i}$ commutes with the Hodge star $*$
 for this metric. In particular, they also commute with the
 $\overline{\6}$-Laplace operator and ${\cal
   L}_{X_i}\alpha$ is a harmonic form for every harmonic
 $\alpha$. Since $X_i^{(0,1)}$ is a complex vector field
 which preserves the structure $I$ and transforms
 $(0,1)$-form into $(0,1)$-form, for a
 $\overline{\6}$-harmonic form $\alpha$, we have ${\cal
   L}_{X_i^{(0.1)}}\alpha^{(0,1)} = \overline{\6} f,$ for
 the function $f = \alpha^{(0,1)}\left(X_i^{(0,1)}\right)$. Since
 ${\cal L}_{X_i^{(0.1)}}\alpha^{(0,1)}$ is harmonic, it
 vanishes. Since we can use any $(0,1)$-vector field
 generated by the action and any harmonic $(0,1)$-form, in
 particular $\sqrt{-1}\alpha,$ we see that the vector fields $X_i$
 preserve the harmonic $(0,1)$-forms. Then, any such form
 has a representation $$\alpha = A_1(\theta_1-i\theta_2) +
 A_2(\theta_3-i\theta_4) + \pi^*(\phi),$$ where $A_i$ are
 pull-backs of functions on the base and $\phi$ is a
 harmonic form on the base $X$. Since $X$ is K3 surface,
 $\phi = 0$. Then, from $d\theta_i = \alpha_i,$ we have
 $\overline{\6}(\theta_1-i\theta_2)=\alpha_1-i\alpha_2$ if
 $\alpha_1-i\alpha_2$ is $(2,0)$-form and 0 if its is
 $(1,1)$. On the other side,
 $\overline{\6}(\theta_3-i\theta_4)=0$, since the other
 characteristic classes are (1,1). As a result, we see that
 $h^{0,1}_{\cal I}(M)=2$, if all curvature forms are
 $(1,1)$ (or instantons) and $h^{0,1}_{\cal I}(M)=1,$ if
 we have one of these forms to be of type $(2,0)$. By
 \ref{_main_intro_Theorem_}, in the first case there is an HKT-metric
 and in the second there is none.

In the construction above we can use a flat 4-tori as a
base instead of $K3$ surface. Then $M$ is a nilmanifold
which corresponds to an example which appeared in
\cite{_Fino_Gra_}. Consider the nilpotent Lie algebra
$\mathbb{R}\times {\mathfrak h}_7,$ where ${\mathfrak h}_7$
is the algebra of the quaternionic Heisenberg group
$H_7$. Its is spanned by the left-invariant vector fields
$e_1,\cdots,e_8$ and is defined by the following relation
on the basis of the dual 1-forms:

$$
\begin{array}{l}
d e ^i  = 0, i = 1, \ldots, 5\\
d e^6 =  e^1 \wedge e^2 +  e^3 \wedge e^4,\\
d e^7 =  e^1 \wedge e^3  - e^2 \wedge e^4,\\
d e^8 =  e^1 \wedge e^4 + e^2 \wedge e^3
\end{array}
$$

On a compact quotient $M=\mathbb{R}\times H_7/\Gamma,$ consider the family \cite{_Fino_Gra_} of complex structures defined via:

$$
\begin{array}{l}
I_t(e^1) = \frac{t-1}{t}e^2, I_t (e^3) = e^4, J_t (e^5) = \frac{1}{t} e^6, J_t (e^7) = e^8,\\
J_t (e^1) = \frac{t-1}{t}e^3, J_t (e^2) =  - e^4, J_t (e^5) = \frac{1}{t} e^7, J_t (e^6)
= - e^8,
\end{array}
$$
for $t\in (0,1)$. Then, for each $t$, $I_tJ_t=-J_tI_t =
K_t$ defines a hypercomplex structure on $M$. Using
averaging argument in \cite{_Fino_Gra_},
it was shown that for $t=\frac{1}{2}$ the structure is HKT
and for $t\neq \frac{1}{2}$ there is no HKT-metric. Here
we provide a different proof using \ref{_HL_for_HKT_Theorem_} and
\ref{_main_intro_Theorem_}. The manifold $M$ has a projection on
$X=T^4$ which makes it a principal bundle with fiber
4-tori and base 4-tori. Then the forms $e^1,e^2,e^3,e^4$
are pull-backs from forms on the base $X$ and the forms
$e^5,e^6,e^7,e^8$ are connection forms in this bundle. So
(up to a constant), the characteristic classes of the
bundle are $0, e^1 \wedge e^2 +  e^3 \wedge e^4, e^1
\wedge e^3  - e^2 \wedge e^4, e^1 \wedge e^4 + e^2 \wedge
e^3$. Now, we note that
\begin{eqnarray*}
d (e^7 + i e^8)&=& (e^1+ie^2)\wedge (e^3+i e^4)\\
&=&\frac{2t-1}{2t-2}\left(e^1+ i\frac{t-1}{t} e^4\right)\wedge(e^3+ i
e^4)\\ && - \frac{1}{2t-2}\left(e^1 - i\frac{t-1}{t} e^2\right)\wedge (e^3 + i e^4).
\end{eqnarray*}
So, when $t=\frac{1}{2}$, it is of type $(1,1)$ with respect
to $I_{\frac{1}{2}}$, but for $t\neq \frac{1}{2}$, it is
of type $(2,0)+(1,1)$. Moreover, the $(2,0)$ component in
this case is  $ \partial_t(e^7 - i e^8) =
\frac{2t-1}{2t-2}(e^1+ i\frac{t-1}{t} e^4)\wedge(e^3+ i
e^4)$, which defines a positive (2,0)-current. So, there
is no HKT-structure if $t\neq \frac{1}{2}$ by
\ref{_HL_for_HKT_Theorem_}. Similarly, we can calculate the Hodge number
$h^{0,1}(M,I_t)$ to check its parity. Instead of using
the fibration structure, its easier to use the result of Console and Fino
(\cite{_Console_Fino_}) who proved that the Dolbeaut cohomology  of a
nilmanifold with an invariant complex structure are
isomorphic to the $\overline{\6}$-cohomology of the
complex of invariant forms. From the defining equations
above, we see that $e^1+ie^2,$ $e^3+ie^4$ and $e^5-ie^6$ are
nonzero elements of $H^{0,1}(M,I_t)$. Also,
$\overline{\6}_t(e^7-ie^8) =  d(e^7-ie^8)|^{(0,2)} =
\frac{2t-1}{2t-2}(e^1- i\frac{t-1}{t} e^4)\wedge(e^3- i
e^4)$. So, for $t=\frac{1}{2},$ it is non-zero in the cohomology and
$h^{0,1}=4$. When $t\neq\frac{1}{2},$ it is not
$\overline{\6}$-closed, $h^{0,1}(M)=3$ and we can apply
\ref{_HL_for_HKT_Theorem_}.

\hfill

{\bf Acknowledgements:}
We are grateful to Dan Popovici for interesting discussions
about strong Gauduchon metrics.

\hfill

\hfill

\hfill

{\small

\noindent
{\sc Gueo Grantcharov\\
{\sc Department of Mathematics and Statistics\\
Florida International University\\
Miami Florida, 33199, USA}\\
\tt grantchg@fiu.edu}\\

\noindent
{\sc Mehdi Lejmi}\\
{\sc D\'epartement de Math\'ematiques\\
\sc Universit\'e Libre de Bruxelles CP218, \\
Boulevard du Triomphe, Bruxelles 1050, Belgique. \\
{\tt mlejmi@ulb.ac.be}\\

\noindent
{\sc Misha Verbitsky}\\
{\sc  Laboratory of Algebraic Geometry, SU-HSE,\\
7 Vavilova Str. Moscow, Russia, 117312}\\
{\tt  verbit@mccme.ru\\}

 }}

\end{document}